\documentclass{amsart}
\usepackage{amssymb}
\newcommand{\mapdown}{\downarrow }
\newcommand{\cal}{\mathcal}
\newcommand{\lo}{\longrightarrow}

\newcommand{\pj}{projection }
\newcommand{\pjs}{projections }
\newcommand{\calg}{$C^{*}$-algebra }
\newcommand{\calgs}{$C^{*}$-algebras }
\newcommand{\walg}{$W^{*}$-algebra }
\newcommand{\walgs}{$W^{*}$-algebras }
\newcommand{\kac}{Kac algebra }
\newcommand{\kacs}{Kac algebras }
\newcommand{\K}{$K=(M,\delta,\kappa,\phi)$ }
\newcommand{\ts}{topological space }
\newcommand{\cpt}{compact }

\newcommand{\hau}{Hausdorff }

\newcommand{\lcpt}{locally compact }
\newcommand{\cnt}{continuous }
\newcommand{\ai}{approximate identity }
\newcommand{\fn}{function }
\newcommand{\fns}{functions }
\newcommand{\homo}{homomorphism }
\newcommand{\shomo}{$*$-homomorphism }

\newcommand{\repn}{representation }
\newcommand{\repns}{representations }
\newcommand{\pro}{projective limit }

\newcommand{\Mpro}{$M=\varprojlim_\alpha M_\alpha$ }

\newcommand{\proc}{pro-$C^*$-algebra }
\newcommand{\procs}{pro-$C^*$-algebras }

\newcommand{\pros}{pro-$C^*$-algebras }
\newcommand{\prow}{pro-$W^*$-algebra }
\newcommand{\prows}{pro-$W^*$-algebras }

\newtheorem{defi}{Definition}[section]
\newtheorem{prop}{Proposition}[section]
\newtheorem{theo}{Theorem}[section]
\newtheorem{lemm}{Lemma}[section]
\newtheorem{cor}{Corollary}[section]
\newtheorem{ex}{Example}[section]

\newtheorem{nota}{Notation}[section]
\begin{document}

\title[local structure of operator
algebras]{non commutative topology and local structure of operator
algebras}
\author[M. Amini]{Massoud Amini}
\address{Department of Mathematics\\University of Shahid Beheshti\\
Evin, Tehran 19839\\Iran\\m-amini@cc.sbu.ac.ir
\linebreak
\indent Department of Mathematics and Statistics\\ University of Saskatchewan
\\106 Wiggins Road, Saskatoon\\ Saskatchewan, Canada S7N 5E6\\mamini@math.usask.ca}
\thanks{* Supported by an IPM grant.}
\maketitle

\begin{abstract}
{\small Starting with a \walg $M$ we use the inverse system obtained by
cutting down $M$ by its (central) projections to define an inverse
limit of \walgs , and show that how this \prow encodes the local
structure of $M$. For the \calgs we do the same thing using their
atomic enveloping \walgs . We investigate the relation of these
ideas to the Akemann-Giles-Kummer non commutative topology.
Finally we use these ideas to look at the local structure of Kac
algebras.}
\end{abstract}

\section{Introduction}

The main motivation of this work is a hope to use the {\it
projective limits} as a device to explain some of the features of
the local theory of non compact quantum groups. However, since
the theory of locally
compact quantum groups is in its beginning [KV], one could not test the applicability of
these ideas quite clearly. Kac algebras were understood as an
appropriate model to generalize the duality of \lcpt topological
groups inside one category. They also provide a solid framework
for quantum groups. They have, however, two shortcomings. One is
the fact that in this theory one assumes (and not {\it prove}) the
existence of the {\it Haar measure}. The other is that we don't
know enough {\it nontrivial} examples of \kacs . 
But they have some advantages also, the most important one being the 
fact that they have a well established theory. In this paper we
use the projective limits of \walgs to study the {\it local}
structure of \kacs . We start with doing the same thing for an
arbitrary \walg and then to adapt the machinery to \calgs. In
particular we study the topological algebra of the elements which
{\it locally} belong to the underlying \walg. We also use the
Pedersen's ideal and its multiplier algebra to make connections
with the already well investigated theory of non commutative
topology.

\section{Pro-$W^*$-algebras}

In the beginning of this section we follow [Frg] to give a
general overview of the {\it
projective limits} of \walgs (For \procs see [Ph88 a-c] and references there in). 
Then we prove some new results which
are needed in the next section. $M$ is called a \prow if \Mpro is an
inverse (projective) limit of \walgs . We denote the
corresponding morphisms of the inverse system
with $\pi_{\alpha\beta} :M_\beta\to M_\alpha$, and $\pi_\alpha : M\to
M_\alpha$ which satisfy
$\pi_{\alpha\beta}\pi_{\beta\gamma}=\pi_{\alpha\gamma}$ and
$\pi_{\alpha\beta}\pi_{\beta}=\pi_{\alpha}$ for each
$\alpha\leq\beta\leq\gamma$.

If \Mpro is a \prow , then it is a \pro and $M=\varprojlim_{\alpha}
\pi_\alpha (M)$. Let's denote the norm and weak operator
topologies of $M_\alpha$ with $\tau_\alpha$ and $\sigma_\alpha$,
respectively. Then $M$ has the {\it projective topology} $\tau
=\varprojlim_{\alpha} \tau_\alpha$. Also if $\pi_{\alpha\beta}$'s are
normal (equivalently $weak^*$-\cnt),
then we can endow $M$ with the $\sigma$-{\it weak projective
topology}  $\sigma=\varprojlim_{\alpha} \sigma_\alpha$ which is {\it
coarser} than $\tau$. In this case we have
$$
(M,\sigma)=\varprojlim_{\alpha}
cl_{\sigma_\alpha} (\pi_\alpha(M))=\varprojlim_{\alpha}
(M_\alpha ,\sigma_\alpha).
$$
Also each $\sigma$-closed $*$-subalgebra of a \prow is a \prow. We are
mostly interested in the case where the morphisms $\pi_{\alpha\beta}$'s are
normal and surjective.

Three immediate examples of \prows are as follows. Let $(M_\alpha)$ be a
{\it descending} chain of \walgs such that the relative norm and weak
topologies are compatible, when restricted appropriately. Then
$$
\cap_\alpha M_\alpha =\varprojlim_{\alpha} M_\alpha
$$
is a (possibly trivial!) \prow . In this case the morphism
$\pi_{\alpha\beta}: M_\beta\to M_\alpha$ is simply the inclusion map
which is normal (by assumption on topologies) but not surjective in
general.

Also for each directed net $(M_\alpha)$ of
\walgs $\prod_\alpha M_\alpha$ is a \prow . Indeed
$$
\prod_\alpha M_\alpha\backsimeq \varprojlim_{\alpha}
(\sum_{\beta\leq\alpha} \bigoplus M_\beta) .
$$
In this case $\pi_{\alpha\beta}:\sum_{\lambda\leq\beta}^{\oplus}
M_\lambda\to \sum_{\lambda\leq\alpha}^{\oplus} M_\lambda$ is just the
projection onto a direct summand which is obviously normal.

The third example is the
most important, since it provides a {\it universal} \prow . Let
$(H_\alpha)$ be a {\it ascending} chain of Hilbert spaces
(i.e. $H_\alpha\subseteq H_\beta$ and $<,>_\alpha =<,>_\beta$ on
$H_\alpha$ for $\alpha\leq\beta$). Then $H=\varinjlim_\alpha H_\alpha
$ with the canonical {\it inductive topology} is called a {\it locally
Hilbert space} [Ino]. Note that $H$ is {\it not} a Hilbert space in
general (it is a Hilbert space if there is a maximal element for the
chain).  Let $\iota_{\alpha\beta}: H_\alpha\to H_\beta$ be
the canonical injection, and define
$$
\cal L (H)=\{T\in L(H): T_\beta\circ \iota_{\alpha\beta}=
\iota_{\alpha\beta} T_\alpha\}
$$
where $T_\alpha\in B(H_\alpha)$ is the restriction of $T$ to $H_\alpha$.
Then $ \cal L (H)=\varprojlim_{\alpha} B(H_\alpha)$ is a \prow and the
morphisms $\pi_{\alpha\beta}: \cal B (H_\beta)\to \cal B (H_\alpha)$ are
normal surjections. We endow $\cal L (H)$ with the $\sigma$-weak
projective topology
$\sigma=\varprojlim_{\alpha} \sigma_\alpha$, where
$\sigma_\alpha=\sigma(\cal B (H_\alpha), K(H_\alpha)^* )$. As an example, if
$H_n=\Bbb C ^n\ \ (n\geq 1)$, then $\cal L (H)$ is just the algebra of
all infinite upper triangular matrices.

The \prow $\cal L (H)$ could be used to give a  {\it \repn theorem}
for \prow which is analogue to the corresponding result for \walgs
. This result asserts that every \prow with the $\sigma$-weak projective
topology is (isomorphic, as topological algebras, to) a
$\sigma$-weakly closed $*$-subalgebra of $\cal L (H)$, for some
locally Hilbert space $H$. Indeed, if \Mpro is a \prow with faithful
weak operator \cnt
representations $\phi_\alpha : M_\alpha \hookrightarrow \cal B (K_\alpha)$,
then we put $H_\alpha=\bigoplus_{\beta\leq\alpha} K_\beta$ and
$H=\varinjlim_\alpha H_\alpha$. Then there is an obvious algebraic
monomorphism $\phi:M\hookrightarrow \cal L (H)$. One can show that the
mappings $\pi_\alpha \phi$ are \cnt , where $\pi_\alpha:\cal L (H)\to
\cal B (H_\alpha)$ are the canonical epimorphisms. Also $\phi(M)\subseteq
\cal L (H)$ is $\sigma$-weakly closed [Frg]. This result
in particular implies that the {\it center}
$$
Z(M,\sigma)= \varprojlim_{\alpha} Z(cl_{\sigma_\alpha} (\pi_\alpha(M)))
$$
is also a \prow .

There is an analogue of the weak operator topology on $\cal L (H)$,
which we denote with $\omega = \varprojlim_{\alpha} \omega_\alpha$,
where $\omega_\alpha=\sigma(\cal B (H_\alpha) ,\cal F (H_\alpha))$. As
$\omega_\alpha \leq \sigma_\alpha$, for each $\alpha$, we have
$\omega\leq\sigma$. In particular, each $\omega$-closed $*$-subalgebra
of $\cal L (H)$ is a \prow . Indeed one can define the analogues of
the {\it \cpt } and {\it finite rank} operators on a locally Hilbert
space in such a way that
$$
\sigma=\sigma(\cal L (H), \cal K (H)) \ \ \text{and} \ \ \omega=\sigma
(\cal L (H), \cal F (H)).
$$

Next we mention some elementary lemmas which are used in the next
section. The proofs are easy and are omitted. Following our
terminology, we use the name {\it pro-Banach space} for the
projective limits of Banach spaces. These results hold in an algebraic
context for vector spaces and pro-vector spaces.

\label{1}
\begin{lemm} Let $E=\varinjlim_{\alpha} (E_\alpha,
\iota_{\alpha\beta})$ be an inductive limit of Banach spaces with
the inductive topology $\tau= \varinjlim_{\alpha} \|.\|_\alpha$.
Then the dual spaces form an inverse system and $$
\varprojlim_{\alpha} E_\alpha ^* \backsimeq(\varinjlim_{\alpha}
E_\alpha)^* , $$ as vector spaces. We denote the pro-Banach space
of the right hand side with $E^*$ and call it the dual of $E$. The
$weak^*$-topology induced from the right hand side on $E^*$ is
denoted by $\tau ^*$. \qed
\end{lemm}

\begin{lemm} Let $( E_\alpha ,\pi_{\alpha\beta})_{\alpha
,\beta \in \Lambda }$ be an inverse system of Banach spaces and
$E=\varprojlim_{\alpha} E_\alpha$,\ $\pi_\alpha : E\to E_\alpha$
be the corresponding pro-Banach space and morphisms. Let
$F_\alpha$ be a closed subspace of $E_\alpha$ such that
$\pi_{\alpha\beta}( F_\alpha)\subseteq F_\beta$, for each $\alpha$
and $\beta$. Then $F=\{x\in E : \pi_\alpha (x)\in F_\alpha \ \
(\alpha\in \Lambda )\}$ is a pro-Banach space and
$F=\varprojlim_{\alpha} F_\alpha$ with morphisms being the
restriction of the $\pi_{\alpha\beta}$'s. \qed
\end{lemm}

\begin{lemm} Consider the inductive limits of Banach spaces 
$E=\varprojlim_{\alpha} (E_\alpha
,e_{\alpha\beta}, e_\alpha)$, $F=\varprojlim_{\alpha} (F_\alpha
,f_{\alpha\beta}, f_\alpha)$, $G=\varprojlim_{\alpha} (G_\alpha
,g_{\alpha\beta}, g_\alpha)$, and $H=\varprojlim_{\alpha}
(H_\alpha ,h_{\alpha\beta}, h_\alpha)$, indexed with the same directed set. 
Then: 
\begin{enumerate}
\item If $\phi_\alpha : E_\alpha \to F_\alpha $ satisfies
$f_{\alpha\beta}\phi_\alpha=\phi_\beta e_{\alpha\beta}$ for each
$\alpha$ and $\beta$, then we can define a map
$\phi=\varprojlim_{\alpha} \phi_\alpha :E\to F$ by $ f_\alpha
\phi=\phi_\alpha e_\alpha$.
\item If $\psi=\varprojlim_{\alpha} \psi_\alpha :F\to G$ is defined
similarly, then $\psi\circ\phi=\varprojlim_{\alpha} \psi_\alpha
\circ\phi_{\alpha}$ is also well defined.
\item If $F\subseteq E$ and $G\subseteq H$ as in previous lemma, and
$\phi=\varprojlim_{\alpha} \phi_\alpha :E\to H$ be as above, then
$\phi(F)\subseteq G$ iff $\phi_{\alpha}(F_{\alpha})\subseteq
G_{\alpha}$ for each $\alpha$.\qed
\end{enumerate}
\end{lemm}
\begin{defi} Given a net $I$, a subnet $J\subseteq I$ is called a
cofinal subnet if for each $x\in I$, there is $y\in J$ with $x\leq y$.
\end{defi}

\begin{lemm} Let $( E_\alpha
,\pi_{\alpha\beta})_{\alpha ,\beta \in I }$ be an inverse system
of Banach spaces (algebras) and $J\subseteq I$ be a subnet of $I$.
Then $( E_\alpha ,\pi_{\alpha\beta})_{\alpha ,\beta \in J }$ is
also an inverse system. Let $E_I=\varprojlim_{\alpha\in I}
E_\alpha$, and $E_J=\varprojlim_{\alpha\in J} E_\alpha$ be the
corresponding pro-Banach spaces (algebras), and consider the well
defined morphism $\pi_{IJ}: E_I\to E_J$ defined by
$(x_\alpha)_{\alpha\in I}\mapsto (x_\alpha)_{\alpha\in J}$. Then
$\pi_{IJ}$ is injective iff $J$ is a cofinal subnet of $I$. If
this is the case, $E_I$ and $E_J$ are isomorphic (through
$\pi_{IJ}$) as topological vector spaces (algebras).
\end{lemm}

To have a coherent notation, we denote each element $x\in
E=\varprojlim_{\alpha} (E_\alpha ,e_{\alpha\beta}, e_\alpha)$ with
$x=\varprojlim_{\alpha} x_\alpha$, where $x_\alpha=e_\alpha (x)$ for
each $\alpha$. We use this notation extensively in the next section.

Now let \Mpro be a \prow for which the morphisms
$\pi_{\alpha\beta}:M_\beta\to M_\alpha$ are normal.
We want to introduce the concept of {\it
predual} for $M $. For $\alpha\leq\beta $, the epimorphism
 $\pi_{\alpha\beta}:M_\beta\to M_\alpha$ induces a monomorphism
$\pi_{\alpha\beta_*}:M_{\alpha_*}\to M_{\beta_*}$. Clearly
$(M_{\alpha_*},\pi_{\alpha\beta_*})$ is an
inductive system. Let $M_* =\varinjlim_\alpha M_{\alpha_*}$. We denote
the corresponding inductive (projective, respectively) topology of $M_*$ ($M$, 
respectively) by
$\tau_*$ ($\tau$, respectively). Then,
by Lemma 2.1, we have:

\begin{prop} $(\tilde M_* ,\tau_*)^* =(\tilde M ,\tau)$.
\qed
\end{prop}

To introduce the concept of the {\it commutant} for $M$, we use
the above mentioned \repn theorem to restrict ourself, without
loss of generality, to the case that $M=\varprojlim_\alpha
(M_\alpha ,\pi_{\alpha\beta})$ is a sub-\prow of $\cal L
(H)=\varprojlim_\alpha (B(H_\alpha) ,\tilde\pi_{\alpha\beta})$ for
some locally Hilbert space $H=\varinjlim_\alpha H_\alpha$, where
$M_\alpha\subseteq B(H_\alpha)$ and the normal morphism
$\pi_{\alpha\beta}$ is the restriction of
$\tilde\pi_{\alpha\beta}$ to $M_\alpha$, for each $\alpha$ and
$\beta$. Then we put $$ M^{'}=\{T\in \cal L (H): Tx=xT \ \ (x\in
M)\}. $$

\begin{prop} $M^{'}$ is a \prow . Moreover, if the
mappings $\pi_{\alpha\beta}$ are surjective, then $M^{'}\backsimeq
\varprojlim_\alpha M_{\alpha} ^{'}$ as \prows .
\end{prop}
{\bf Proof} $M^{'}$ is obviously a $\sigma$-closed $*$-subalgebra
of $\cal L (H)$, so it is a \prow . Next assume that
$\pi_{\alpha\beta}$'s are surjective. Then
$\tilde\pi_{\alpha\beta}( M_{\beta} ^{'})\subseteq M_{\alpha}
^{'}$ (for each $x\in M_\beta$ and $y\in M_\beta ^{'}$, we have $
\tilde\pi_{\alpha\beta}(y)\pi_{\alpha\beta}(x)
=\tilde\pi_{\alpha\beta}(yx)=\tilde\pi_{\alpha\beta}(xy)
=\pi_{\alpha\beta}(x)\tilde\pi_{\alpha\beta}(y)$, and each element
of $M_\alpha$ is of the form $\pi_{\alpha\beta}(x)$ for some $x\in
M_\beta$, so $ \tilde\pi_{\alpha\beta}(y)\in M_{\alpha} ^{'}$).
Let $\pi_{\alpha\beta} ^{'}$ to be the restriction of
$\tilde\pi_{\alpha\beta}$ to $M_{\alpha} ^{'}$, then $(M_\alpha
^{'} ,\pi_{\alpha\beta} ^{'})$ is an inverse system and we can
form the inverse limit as a sub-\prow $\cal L (H)$. But
$\tilde\pi_{\alpha}(M^{'})\subseteq \{T_\alpha\in B (H_\alpha):
T_\alpha x_\alpha=x_\alpha T_\alpha\ \ (x_\alpha\in M_\alpha)\}$,
which is clearly contained in $M_{\alpha} ^{'}$. Hence
$M^{'}\subseteq \{T\in \cal L(H): \tilde\pi_\alpha(T) \in M_\alpha
^{'}\ \text{for each} \alpha\}$. Conversely, if $T$ belongs to the
RHS, then given $x\in M$, $ \tilde\pi_\alpha(Tx-xT)=0$, for each
$\alpha$, and so $Tx-xT=0$, which means that $T\in M^{'}$. Now by
Lemma 2.3, we get $M^{'}\backsimeq \varprojlim_\alpha
M_{\alpha} ^{'}$ as sub-\prows of $\cal L (H)$. \qed

It is natural to ask if a generalization of {\it von Neumann double
commutant theorem} holds for \prows . The answer is affirmative.

\begin{theo}{\bf (Generalized double commutant theorem)}
Assume that $M=\varprojlim_\alpha (M_\alpha , \pi_{\alpha\beta})$ is a
\prow such that the morphisms $\pi_{\alpha\beta}$ and
$\pi_{\alpha\beta}^{'}$ are surjective, for each $\alpha$ and
$\beta$. Then $M^{''}$ is a \prow and $M\backsimeq M^{''}$ as
\prows .
\end{theo}
{\bf Proof} By applying the above proposition twice, we have
$M^{'}\backsimeq \varprojlim_\alpha M_{\alpha} ^{'}$ and
$M^{''}\backsimeq \varprojlim_\alpha M_{\alpha}
^{''}=\varprojlim_\alpha M_{\alpha}=M$. \qed

\section{Non commutative topology}

It is the basic philosophy of the non commutative topology to
think of {\it projections} as analogues of {\it sets}. Therefore
one need to use \walgs , where it is guaranteed that we have
plenty of projections. Along this line of thinking, one would
consider the {\it minimal projections} as {\it points}. To have
enough points, then one restricts to the case of an {\it atomic
\walg} in which every projection majorizes a minimal one. Yet, not
to loose connection with {\it topological structure}, we have to
have this \walg comming from a $C^*$-algebra. There are two
(basically equivalent) way of proforming this task. We start with
a \calg $A_0$ and consider the {\it atomic \repn} $\{\pi_a ,H_a\}$
of it, namely the direct sum of a maximal family $\{\pi_\alpha
,H_\alpha\}$ of mutually non equivalent irreducible \repns of
$A_0$ (not unique, but spatially equivalent with any other
choice), and put $A_a =\pi_a(A_0)^{''}\subseteq B(H_a)$. This is
an atomic \walg and (a copy of) $A_0$ is strongly dense in it
[Dix]. A more accessible approach, which we adopt in this section,
is to consider the union $z$ of all minimal projections in
$A_0^{**}$, then $z$ is a central projection and $A_a =z A_0^{**}$
is an atomic \walg which contains a strongly dense copy of $A_0$,
namely $zA_0$ [Ak69]. If $X$ is a \lcpt \ts and $A_0 =C_0(X)$,
then $A_a=B(X)$ is the algebra of all bounded functions on $X$ (if
$\mu$ is a Borel measure on $X$, then $L^\infty (X)$ is the
quotient of $B(X)$ obtained by identifying \fns which are equal
$\mu$-a.e.).

The non commutative topology is then a certain subset of the lattice
$Prj(A_a)$ of projections of $A_a$. Before going into this, let's
introduce some notations due to G. Pedersen. For $E\subseteq A_a$, let
$E^m$ be the set of {\it suprema} in $(A_a)_{sa}$ (self adjoint part
of $A_a$) of all norm bounded increasing nets of elements of
$E_{sa}$. Also let $E_m=-(-E)^m$. We say that $E$ is {\it lower
monotone closed} if $E^m=E$. Now
$$
\tau=Prj_{op}(A_a)=Prj(A_a)\cap (A_0 ^+ )^m
$$
is called the {\it Akemann-Giles-Kummer topology} and its elements
are called {\it open} projections. The complements of open
projections are called {\it closed} and their collection is denoted by
$Prj_{c\ell}(A_a)$. For each \pj $p$, its {\it closer} ({\it interior}
) is defined to be the infimum (supremum) of all closed (open)
projections majorizing (majorized by) $p$, which is a closed (an open)
\pj and is denoted by $\bar p$ ($p^\circ$). Also the {\it central
cover (support)} of $p$ is denoted by $z(p)$ and is defined to be the
infimum of all central \pjs majorizing $p$.
We say that a \pj $p\in Prj(A_a)$ is {\it
(relatively) compact} if it is closed and $pa=p$ (or equivalently
$p\leq a$) for
some $a\in A_0 ^+$ (if $\bar p$ is \cpt, resp.), and
denote their collection by $Prj_{cp}(A_a)$ ($Prj_{rc}(A_a)$, resp.).
$p\in Prj(A_a)$ is called
{\it regular} if $\|pe\|=\|\bar p e\|$, for each $e\in Prj_{op}(A_a)$
, and {\it quasi-compact}, if for each net $(b_\alpha)\subseteq
U(\tau)^+$ (see next paragraph for a definition!) with $b=\wedge
b_\alpha \in U(\tau)^+$, $inf_\alpha \|pb_\alpha p\|=\|pbp\|$. We use
the notations $Prj_{rg}(A_a)$ and $Prj_{qc}(A_a)$ to denote the
collection of these projections. Also we denote the collection of all
central and minimal projections by $Prj_{cn}(A_a)$ and
$Prj_{mi}(A_a)$, respectively. The last set should be understood as
{\it points} of our space.

One should note that \cpt \pjs may not have the {\it finite covering
property }, namely if $p$ is a \cpt \pj and $p\leq \bigvee p_i$, where
$p_i$'s are open projections, one may not be able to {\it cover} $p$
with a finite subfamily of $p_i$'s. A contraexample could simply be
constructed when $p$ and $p_i$'s are rank one \pjs on a Hilbert space
[Ake]. However, a version of {\it finite intersection property} holds
for \cpt \pjs : If $p$ is \cpt and $p_i$'s are closed \pjs such that
any finite intersection of \pjs $p\wedge p_i$ is non zero, then
$p\wedge (\wedge_i p_i)\neq 0$.

\begin{nota} Most of the time we need to specify
that a \pj belongs to two or three classes. We suggest the
following abbreviations to avoid lengthy indices. To show a single
class of operators we use two letters from the corresponding
class's name, so as before: cn, op, c$\ell$, cp, rc, qc, and rg
stand for central, open, closed, compact, relatively compact,
quasi compact, and regular, respectively. Now we use also: gq, ro,
co, cm, and cro to indicate the class of \pjs which are "regular
and quasi compact", "relatively compact and open", "central and
open", "central and minimal", and "central, relatively \cpt , and
open", respectively. For instance, $Prj_{ro}(M)=Prj_{rc}(M)\cap
Prj_{op}(M)$.
\end{nota}

With this notation we have
$Prj_{gq}(A_a)=Prj_{cp}(A_a)$ [Sha, 5.3], and
$Prj_{cp}(A_a)\subseteq Prj_{c\ell}(A_a)$[Sha, 5.2].

A self adjoint element $a\in A_a$ is called {\it continuous} (with
respect to $\tau$) if the {\it spectral projections} $E_a(U)$ of
$a$ corresponding to open sets $U\subseteq \Bbb R$ are open
projections. It is said to {\it vanish at infinity} if, moreover,
$E_a(K)$ is a \cpt \pj for each compact set $K\subseteq \Bbb R$.
The set of all elements of $A_a$ whose real and imaginary parts
are \cnt (and vanishing at infinity) is denoted by $C_b(\tau)$
($C_0(\tau)$, respectively). Let $\Lambda^+ (\tau)$ be the minimal
convex cone containing $\tau$ and $\Lambda (\tau)=\Lambda^+
(\tau)+\Bbb R 1$. Then the minimal lower monotone closed subset of
$(A_a)_{sa}$ containing $\Lambda(\tau)$ is denoted by $L(\tau)$
and its elements are the {\it lower semi \cnt} (self adjoint)
elements of $A_a$. The {\it upper semi \cnt} (self adjoint)
elements of $A_a$ are defined similarly and denoted by $U(\tau)$.
An element of $A_a$ whose real and imaginary parts are both lower
and upper semi \cnt is called {\it quasi \cnt} and $Q(\tau)$
denotes their collection. Then we have the following equalities:
$$ A_0=C_0(\tau),\ \ A_b=M(A_0)=A_0 ^m\cap A_{0_m}=C_b(\tau),\ \
\ \ QM(A_0)=\overline{\tilde A_0 ^m}\cap
\overline{\tilde A_{0_m}}=Q(\tau), $$ where $QM(A_0)$ is the space
of quasi multipliers of $A_0$ [AR].

The Akemann-Giles-Kummer (quasi) topology $\tau$ on the atomic \walg $M=A_a$
satisfies a number of axioms (it is so called a $C^*$-{\it topology}):
\begin{enumerate}
\item $0,1\in \tau$,
\item For each $(p_\alpha)\subseteq \tau \ \ \vee p_\alpha \in \tau$,
\item If $p,q\in \tau$ and $pq=qp$ then $p\wedge q\in \tau$,
\item $\tau=L(\tau)\cap Pr(M)$,
\item For each $a\in \Lambda^+ (\tau)$, \ $a^{\frac1 2}\in L(\tau)$,
\item For each unitary $u\in C_b(\tau)$, \ $u^* \tau u\subseteq \tau$.

Moreover it has the following properties:
\item $T_1$: $\forall e,f\in Prj_{mi}(A_a)\ \exists p\in Prj_{op}(A_a)\
\ \ ef=0\ \vdash \ e\leq p\ \text{and}\  pf=0$,
\item {\it Hausdorff}: $\forall e,f\in Prj_{mi}(A_a)\ \exists p,q\in
Prj_{op}(A_a)\ \ \ ef=0\ \vdash \ e\leq p\  \text{and}\  f\leq q$,
\item {\it locally \cpt}: $\forall e\in Prj_{mi}(A_a)\ \exists p\in
Prj_{ro}(A_a)\ \ \ e\leq p$,
\item {\it regular}: $\forall e\in Prj_{mi}(A_a)\ \forall f\in
Prj_{c\ell}(A_a)\ \exists p,q\in Prj_{op}(A_a)\
\ \ ef=0\ \vdash \ pq=0, \ p\geq e,\ \text{and}\  q\geq f$ ,
\item {\it completely regular}: $\forall e\in Prj_{mi}(A_a)\ \forall f\in
Prj_{c\ell}(A_a)\ \exists a\in C_0(\tau)_1 ^+
\ \ ef=0\ \vdash \ ae=e\ \text{and} \ af=0$,
\end{enumerate}
One can show that a $C^*$-topology (on any atomic \walg) is
the Akemann-Giles-Kummer topology of a \calg iff it is $T_1$, completely
regular , and \lcpt
[Sha, 5.5]. We say that $\tau$ is {\it \cpt} if $1\in
Prj_{cp}(A_a)$. This holds iff $A_0=A_b=M(A_0)$ (that's when $A_0$ is
unital) [Sha,5.6].

Although $A_a$ provides an ideal setting for {\it
topological} observations, it is still quite {\it large}
for measure theoretic purposes (if $A_0=C_0(X)$ then the \fns in
$A_a=B(X)$ are not necessarily measurable). On the other hand, any
sucessful theory of locally compact (quantum) groups should appreciate
the close relation between {\it topology} and {\it measure
theory}. This could be demonestrated by the fact that for a
\hau topological group (of second category) the existence of the Haar
measure and locally compactness are indeed equivalent [Ng].

Fortunately there is a well established theory to be used for this
purpose, namely the {\it measurable} decomposition theory [Pd79, ch4].
We quote here some of the basics of this theory. Recall that if
$\{\pi_a , H_a\}=\oplus_{\alpha\in \hat A_0} \{\pi_\alpha ,H_\alpha\}$
be the atomic \repn of a \calg $A_0$, then
$A_a=\pi_a(A_0)^{''}\subseteq B(H_a)$ is a atomic \walg .
Indeed $A_a =\sum_{\alpha} ^\oplus B(H_\alpha)$. Now $\pi_a$ has a unique
extension to a normal \repn $\bar\pi_a$ of $A_0^{**}$, such that
$\bar\pi_a (A_0^{**})=\sum_{\alpha} ^\oplus B(H_\alpha)$. Hence, $\pi_a$
being faithful, this shows that $\bar\pi_a$ is not faithful in
general. However, one can single out a large class inside $A^{**}$ on
which  $\bar\pi_a$ is faithful (isometric). An element $x\in
A_{sa}^{**}$ is called {\it universally measurable} if for each
$\epsilon >0$ and each $\omega\in S(A_0)$, there are $a,b\in
((A_0)_{sa})^m$ such that $-b\leq x \leq a$ and
$\omega(a+b)<\epsilon$. The family of all such elements is denoted by
$\cal U (A_0)$. This is a norm closed real vector space containing
$(({\tilde A}_{sa})^m )^{-}$. Thus
$A_{mu}=\cal U (A_0) +i \cal U (A_0)$ is a $C^*$-subalgebra of
$A_0^{**}$. Also elements of $\cal U (A_0)$ could be approximated
strongly from above (below) by elements of $(((A_0)_{sa})^m )^{-}$\
(respectively $(((A_0)_{sa})_m )^{-}$).

Usually it is more useful if one restricts to the class of the
universally measurable ones. For this reason, we consider the real
vector space $\cal B ((A_0)_{sa})$ which is defined as the smallest
class (in $A_0^{**}$) which contains both $(A_0)_{sa}$ and the strong
limits of all monotone (increasing or decreasing) sequences of
elements of $ (A_0)_{sa}$. Then this is the self adjoint part of a
\calg . The \calg $A_{bm} =\cal B ((A_0)_{sa}) +i\cal B ((A_0)_{sa})$
is called the {\it enveloping Borel} $*$-{\it algebra} of $A_0$
(first introduced by Richard Kadison). If
$A_0$ is separable, then $A_{bm}$ is unital (this might fail in
non separable case). In this case, $A_{bm}$ is closed under bounded
Borel functional calculus, central support (cover), and polar
decomposition. Also each \repn
$\{\pi_0, H_0\}$ of $A_0$ uniquely extends to a sequentially normal \repn
$\pi_{bm}$ of $A_{bm}$, and $\pi_{bm} (A_{bm})=(\pi_0(A_0))_{bm}$ in
$B(H_0)$. In particular, when $A_0$ and $\{\pi_0, H_0\}$ are separable, then
$\pi_0^{''}(A_{bm})=\pi_0(A_0)^{''}$ (in the universal Hilbert space),
and center is mapped onto the center.

We know that the atomic \repn is faithful on $A_{um}$. On the
other hand $A_{um}$ is monotone sequentially closed, so in particular
$A_{bm}\subseteq A_{um}$ and $\pi_a$ is also faithful on
$A_{bm}$. Therefore, after appropriate identifications, we have
$$
A_0\subseteq A_{bm}\subseteq A_{um}\subseteq A_{a}\subseteq A_{u}.
$$
In the commutative case, when $A_0=C_0(X)$, $A_{bm}$ and $A_{um}$ are
the algebras of all bounded Borel and universally measurable \fns on
$X$, respectively, where as $A_{a}=A_{u}=B(X)$, the set of all bounded
\fns on $X$.

There are also the strong and weak version (in contrast with the norm
version) of the enveloping Borel $*$-algebra. Consider $A_0$ in its
universal \repn . Then one can define the \calgs $A_{bm}^s$ and
$A_{bm}^w$ similar to $A_{bm}$ by just replacing the norm sequential
monotone limits with strong and weak one (the later is what is called
a $\Sigma^*$-algebra by E.B. Davis [Dav]). Then we have
$$
A_{bm}\subseteq A_{bm}^s\subseteq A_{bm}^w .
$$
It seems that these three \calgs should be indeed equal
(like in the commutative case), but I am not
aware of any general result in this direction (we know an affirmative
answer for the case of type $I$ \calgs ). Now $A_{um}$ is indeed
strongly sequentially closed. Therefore we can toss in these algebras
in the above sequence of inclusions as
$$
A_0\subseteq A_{bm}\subseteq A_{bm}^s\subseteq A_{um}\cap
A_{bm}^w\subseteq A_{um}\cup A_{bm}^w\subseteq A_{a}\subseteq A_{u}.
$$

As for the relation of these algebras with the non commutative
topology, first observe that all open and closed projections of
$A_a$ are limits of monotone nets of (positive) elements of $A_0$.
Therefore if we restrict ourselves to the separable case, we get
$\tau=Prj_{op}(A_a)\subseteq A_{bm}$. Also all closed (and so
compact) projections of $A_a$ already belong to $A_{bm}$.
Therefore one can simply apply most of the {\it topological}
results proven for $A_a$ to $A_{bm}$ also. We make this more
precise in the forthcoming sections. Here we rather want to
emphasize on the importance of $A_{bm}$ in non commutative measure
theory. This is crucial if one wants to deal with unbounded
measures. For instance for a \lcpt group $G$ one can extend the
Haar measure to a weight on $B^\infty(G)$, but it is not possible
(or at least clear how) to do this for $B(G)$.

First let us note that $A_{bm}$ is inside the atomic \walg $A_a$, so
one can hope to have minimal \pjs inside $A_{bm}$. Indeed there are
plenty of them. Let $A_0$ be a separable \calg in its atomic \repn
and $P(A_0)$ be the
set of its pure states (the ones which correspond to irreducible
\repns). Let $\hat A_0$ be the {\it spectrum} of $A_0$, that's the set
of equivalence classes of irreducible \repns of $A_0$. Then there is a
one-one correspondence between $P(A_0)$\ ($\hat A_0$, respectively) and
$Prj_{mi}(A_{bm})$ \ ($Prj_{cm}(A_{bm})$, respectively) [Pd71]. In
this correspondence a state $\omega\in P(A_0)$ is sent to the one
dimensional projection $p$ onto the subspace of $H_a$ spanned by the
unit vector in $H_a$ which represents $\omega$ in its {\it GNS}-\repn
. Conversely given $p\in Prj_{mi}(A_{bm})$, the reduced algebra
$(A_{bm})_p$ is isomorphic to $\Bbb C$ and $\omega$ is just defined on
$A_{bm}$ by $\omega(x)=pxp$, which turn out to be a $\sigma$-normal
pure state of $A_{bm}$, and thereby a pure state of $A_0$. Also the
center of $A_{bm}$ is indeed equal to $B_\sigma ^\infty (\hat A_0)$,
where index $\sigma$ tells us that we have to restrict to a sub
$\sigma$-field of Borel measurable sets (namely the {\it Davis-Borel
structure} on $\hat A_0$), and the second correspondence could be
rephrased as the assertion that the points of $\hat A_0$ are Davis-Borel
sets.

Every state of $A_0$ extends uniquely to a $\sigma$-additive
functional on $A_{bm}$. One can not hope to get a normal state
extension on $A_{bm}$ (neither on $A_a$) in general. For this reason,
G.K. Pedersen has singled out a class of states (which he calls {\it
atomic} states) for which this is possible [Pd71]. Here we briefly quote
this theory. Let's recall that $A_{bm}$ could be embedded in the direct sum
of $B(H_i)$, the bounded operators on the Hilbert space of the \repns
indexed by, and associated with, elements of $P(A_0)$. Let $J$ be the
smallest monotone closed $C^*$-subalgebra of $A_{bm}$ containing
$Prj_{mi}(A_{bm})$. $J$ could be identified with the ideal consisting
of those $(x_i)$'s in the
above direct sum which have countably many nonzero components. A
positive functional $\omega\in A_0 ^{*}$ is called {\it atomic} ({\it
diffuse}) if there is (for all) $p\in Prj_{mi}(A_{bm})$ such that
$\omega(1-p)=0$ ($\omega(p)=0$). Each positive functional then
uniquely decomposes as the sum of an atomic and a diffuse one (which
are centrally orthogonal). The atomic functionals are exactly those
which could be decomposed as a countable linear combination of
(mutually orthogonal) pure states. These functionals extend uniquely
to a normal state on $A_{bm}$ (and $A_a$).

It is useful to note that the above observations have unbounded
analogues. Indeed all of the above statements are also valid for
the following situation [Pd71]: Let $\omega$ be a $\sigma$-normal
weight on $A_{bm}$ which is majorized by an invariant (under
spatial inner automorphisms) convex functional $\rho$ on
$A_{bm}^{+}$ such that there is a sequence $(e_n)$ of
$\rho$-finite elements (\pjs) in $A_{bm}$, summing to 1. These
conditions hold (with $\rho=\omega$), in particular, when $\omega$
is a $\sigma$-finite $\sigma$-trace [Dav], or a $C^*$-integral [Pd
69]. Then $\omega$ decomposes as the countable sum of
$\rho$-normal bounded linear functionals $\omega_n=\omega(e_n .)$.
Moreover, there is a sequence $(p_n)$ of central \pjs in $J$ such
that $\omega_n(p_n .)$\ \ ($\omega_n((1-p_n).)$, respectively) is
atomic (diffuse, respectively). In particular, $\omega$ decomposes
as the sum of the atomic and diffuse weights  $\omega(p.)$ and
$\omega((1-p).)$, where $p=\vee p_n\in Z(A_{bm})\cap J$. As for
the decomposition of an atomic weight $\omega$ into sum (mutually
orthogonal) of pure states, one should note that this is not
doable in general (even with above restrictions) unless $\omega$
is a trace.

\section{The local algebra and the quasi local space of a \calg }

Let $M$ be a \walg on a Hilbert space $H$. If a
bounded operator commutes with elements of
$M^{'}$, then it should belong to $M$ by double commutant
theorem. A closed operator, however, could commute with elements
of $M^{'}$ but does not belong to $M$. In this situation we say that it
is {\it affiliated } with $M$. The set of all closed (densely defined)
operators on $H$ affiliated with $M$ is denoted by $M^\eta$.
There are many situations that we need
to consider affiliated elements. In this section we consider a class
of such elements coming from \prows .

The basic idea is to use {\it cut down} by projections as a
non commutative analog of {\it neighborhoods}. In particular, working
with unbounded elements, we are mainly interested in \cpt
neighborhoods.
There is a notion of compactness of projections defined by
Akemann [Ak69]. This was originally manufactured to get a
Stone-Weierstrass theorem for \calgs , but it also suits our
purposes. We briefly discussed this in previous section. Let's recall
the basic terminology. Let $A_0$ be a \calg and
$A_a=zA^{**}$, where $z\in Prj_{cn}(A^{**})$ is the supremum of all
minimal projections.
In particular $A_a$ is an atomic  $W^*$-algebra. For any closed left ideal
$I\subseteq A_0$ the weak closure of $I$ is a weakly closed left ideal
of $M$ (after identifying $A_0$ with $zA_0$), and so is of the form $A_0p$,
for some $p\in Prj(A_a)$. The projections obtained in this way are
called {\it open} and we denote the set of all open
projections by $Prj_{op}(A_a)$. The complement of these type of
projections are called {\it closed} and are collected in
$Prj_{cl}(A_a)$. Given $p\in Prj(A_a)$, the smallest closed projection
$\geq p$ is called the {\it closure} of $p$ and is denoted by $\bar
p$. A projection $p\in Prj(A_a)$ is open (closed) iff
there is an increasing (decreasing) net of positive elements in $A$
converging to $p$ in $weak^*$-topology. A closed projection $p\in
M$ is called {\it \cpt} if there is $a\in A_0^+$ with $ap=p$. A
projection $p$ is called {\it relatively \cpt} if $\bar p$ is \cpt . We
denote the set of all (relatively) \cpt projections by $Prj_{cp}(A_a)$
($Prj_{rc}(A_a)$, respectively). Given a closed left ideal $I$ of
$A_0$, the corresponding open \pj is relatively compact iff
$I\subseteq A_{00}=K(A_0)$. A closed
projection is compact iff it is closed projection in $\tilde A_0^{**}$,
where $\tilde A_0$ is the minimal unitization of $A_0$ [Ak69]. In particular
all closed projections of $A_0^{**}$ are compact iff $A_0$ is unital.

\begin{lemm} With above notations, TFAE \begin{enumerate}
\item $1\in Prj_{cp}(A_a)$,
\item $Prj_{cp}(A_a)=Prj_{c\ell}(A_a)$,
\item $1\in Prj_{ro}(A_a)$,
\item $Prj_{ro}(A_a)=Prj_{op}(A_a)$,
\item $A_0$ is unital,
\item $A_0=A_0^\eta$,
\item $A_{00}=A_0=A=A_b=A_{q}$.
\end{enumerate}
Moreover, (any of) these conditions would imply that

$(8)$\ \ $H_{q\ell}=H_{\ell}=H_a$, $ \cal L (H_{q\ell})
=\cal L (H_{\ell})=B(H_a)$, and $A_{q\ell}=A_\ell=A_a$.
\end{lemm}
{\bf Proof} The equivalence of $(1)$ and $(2)$ is obvious. Since
$1$ is both open and closed, $(1)$ is equivalent to $(3)$ and
$(4)$ implies $(1)$. Also $(2)$ obviously implies $(4)$. The
equivalence of $(5)$ and $(1)$ is [Sha, 5.6]. Also the equivalence
of $(5)$, $(6)$, and $(7)$ is obvious. Finally $(1)$ implies
$(8)$, because $1$ being a (central) open relatively \cpt \pj, we
get $H_{q\ell}=H_a$ and all the locally bounded operators become
bounded. \qed

For each $p\in Prj_{ro}(A_a)$, $A_p ^a=(A_a)_p=pA_a p$ is a
\walg and we have an
(involutive) linear map $\pi_p ^a:A_a\to A_p ^a,\ x\mapsto pxp$.
Consider the canonical
order on $Prj_{ro}(A_a)$, that's $p\leq q$ iff $pq=qp=p$.
In particular, for $p\leq q$, the
linear map $\pi_{pq} ^a:A_q ^a\to A_p ^a;\ qxq\mapsto pxp$
is well defined and
$$
\pi_{pq}^a\pi_{qr}^a=\pi_{pr}^a,\ \
\pi_{pq}^a\pi_q^a=\pi_p^a\ \ \ (p\leq q\leq r).
$$

Hence we have an inverse system $(A_p ^a, \pi_{pq}^a)$ of \calgs
. Now let
$A_{q\ell}=\varprojlim_{p\in Prj_{ro}(A_a)} A_p^a$ be
the corresponding \proc and denote the
corresponding morphisms from $A_{q\ell}$ to $A_p^a$ by $\pi_p^\ell$, then
$\pi_{pq}^a\pi_q^\ell =\pi_p^\ell$, for each $p,q\in Prj_{ro}(A_a)$.
$A_{q\ell}$ is called the {\it quasi local space} of $A_0$. Note that
this is merely a (locally convex) topological vector space and not an
algebra. Although $A_a^p$'s are all \walgs , but $\pi_{pq}^a$'s are
hardly ever homomorphisms. Indeed even in the trivial case where
$A_0=M_2(\Bbb C), q=\left(\matrix 1&0\\
0&1\endmatrix\right), \text{and}\ p=\left(\matrix 1&0\\
0&0\endmatrix\right)$,
the map $\pi_{pq}^a:x\mapsto pxp$ is not a homomorphism (simply there
is no nonzero algebra homomorphism from $M_2(\Bbb C)$ to $pM_2(\Bbb
C)p$).

Recall that $A_a$ could be considered as $A_a=\pi_a(A_0)^{''}\subseteq
 B(H_a)$, where $\{\pi_a ,H_a\}$ is the {\it atomic \repn} of
$A_0$. Now, for each $p\in Prj_{ro}(A_a)$, consider the closed
subspace $H_p^a=pH_a$ of $H_a$.
If $p,q\in Prj_{ro}(A_a)$ and $p\leq q$, then $H_p^a=p
H_a=qpH_a\subseteq qH_a=H_q^a$. Also the inner product of $
H_q^a$ restricted to $H_p^a$ clearly coincides with that of $
H_p^a$. Hence, if $\iota_{pq}^a: H_p^a \to H_q ^a$ be the
corresponding embedding, we get the local Hilbert space
$H_{q\ell}=\varinjlim_{p\in Prj_{ro}(A_a)} (H_p^a
,\iota_{pq}^a)$ [Ino]. Again note that we are taking the {\it algebraic}
direct limit, so in this case, $H_{ql}$ is indeed the union of
subspaces $pH_a$ of $H_a$, where $p$ runs over all open relatively
\cpt \pjs in $A_a$.
In particular $H_{ql}\subseteq H_a$. The case where this is indeed
dense is more interesting as we see in a moment. But first let's
observe that $A_{ql}$ is a subspace of
$$
\cal L (H_{q\ell}) =\varprojlim_{p\in Prj_{ro}(A_a)} B(H_p).
$$

\begin{prop} With above notations,
\begin{enumerate}
\item For each $p\in Prj_{ro}(A_a)$, $B(H_p^a)=pB(H_a)p$.
\item The linear maps $\pi_{pq}^a:A_q^a\to A_p^a$
are the restrictions of the maps
$\bar\pi_{pq}^a :B(H_q^a)\to B(H_p^a)$ and
$A_{q\ell}$ is a subspace of $\cal L
(H_{q\ell})$.
\item The mapping $x\mapsto (pxp)_{p\in Prj_{cp}(A_a)}$ is an embedding of
$A_a$ into $A_{q\ell}$.
\end{enumerate}
\end{prop}
{\bf Proof} Recall that $H_p^a=pH_a$. Now for each $T\in B(H_a)$
and $\zeta\in H_a$, we have $pTp(p\zeta)=pTp(\zeta)=p(Tp\zeta)$,
so $pTp\in L(H_p^a)$. Also clearly $\|pTp\|_{H_p^a}\le
\|p\|^2\|T\|\leq \infty$, hence $pTp\in B(H_p^a)$. Conversely,
consider the orthogonal decomposition $H_a=pH_a \bigoplus
(1-p)H_a$.  Then $pH_a=\{p\eta :\eta\in pH_a\}$. (The advantage of
choosing $\eta$ in $pH_a$\ - and not in $H_a$\ - is that each
element of $pH_a$ can be written {\it uniquely} in the form
$p\eta$, where $\eta\in pH_a$). Now let $S\in B(H_p^a)$ and
$\zeta\in H_a$, Let $\eta$ be the unique element of $pH_a$ with
$S(p\zeta)=p\eta$. Then the map $T$ which sends $\zeta$ to $\eta$
is a well-defined linear map on $H$ and
$pTp(p\zeta)=pT(p\zeta)=p(p\eta)=p\eta=S(p\zeta)$, for each
$\zeta\in H_a$, i.e. $S=pTp$. To show that $T$ is bounded, observe
that $T(\zeta)=\eta=p\eta=Sp(\zeta)$, so $\|T\|\leq\|S\|_{
H_p^a}\|p\|\le\infty$. Hence $T\in B(H_a)$ and $(1)$ is proved.
Now $(2)$ follows immediately and $(3)$ is trivial.\qed

Now assume for a moment that $\cup_p pH_a$ is norm dense in $H_a$.
Consider $\bar x\in A_a$, and let $x$ be the restriction of
$\bar x$ to $H_a\subseteq H_{q\ell}$. Put
$$
D(x)=\{\zeta\in H_a: x\zeta\in H_a\},
$$
then, for each $p\in Prj_{ro}(A_a)$ and $\zeta\in H_a$, we have
$pxp\zeta\in pH_a \subseteq H_a$, that's $pH_a\subseteq D(x)$. In
particular it follows from the above lemma that $D(x)\subseteq H_a$ is
norm dense. Therefore $\check x: D(x)\subseteq H_a\to H_a$ is a densely defined
(unbounded) operator on $H_a$. Also note that $p\bar x p\zeta=pxp\zeta$,
therefore $\bar x$ is uniquely determined by $x$. Also $pxp=p\check
x p\in A_p^a\subseteq A_a$ for each $p\in Prj_{ro}(A_a)$.

We claim that $\check x\in A_a^\eta$. Indeed,
for each $p\in Prj_{ro}(A_a)$ and $y\in A_a^{'}$, we have
$p(xy-yx)p=y(pxp)-(pxp)y=0$, so $xy-yx=0$ in $\cal L
(H_{q\ell})$, and we are done. Now
identifying $\check x$ with $x$ we have

\begin{prop} If $H_{ql}$ is norm-dense in $H_a$,
then $A_{q\ell}\subseteq A_a^\eta$.\qed
\end{prop}

Up to now we have worked with a variety of projections. It is the time
to see what happens if, in each case, we restrict ourself to the ones
which are also {\it central}. The fact that each projection in a \walg
has a {\it central cover}, and if the \pj is \cpt then so is its
central cover, also motivates us to consider this
situation. In particular we are interested in the
family $Prj_{cro}$ of all
central, relatively \cpt open \pjs of $A_a$.

Now consider the sub-system $(A_p^a,\pi_pq ^a)_{p,q\in
Prj_{cro}(A_a)}$ and form the corresponding \prow $A_\ell
=\varprojlim_{p\in
Prj_{cro}(A_a)} A_p ^a$. This is called the {\it local algebra} of
$A_0$. We have gathered some of the properties of $A_{\ell}$ in
the following proposition which could be proved similar to
the quasi local case.

\begin{prop} Let $A_0$ be a \calg , and $\{\pi_a ,H_a\}$ be its
atomic \repn . Let $A_a=\pi_a(A_0)^{''}$ and $A_\ell$
be the corresponding atomic
\walg and local algebra of $A_0$, then we have
\begin{enumerate}
\item For each $p\in Prj_{cro}(A_a)$, consider the closed
subspace $H_p^a=pH_a$ of $H_a$.
If $p,q\in Prj_{cro}(A_a)$ and $p\leq q$, then $H_p^a=p
H_a=qpH_a\subseteq qH_a=H_q^a$. Also the inner product of $
H_q^a$ restricted to $H_p^a$ coincides with that of $
H_p^a$. If $\iota_{pq}^a: H_p^a \to H_q ^a$ be the
corresponding embedding, then $H_{\ell}=
\varinjlim_{p\in Prj_{cro}(A_a)} (H_p^a
,\iota_{pq}^a)$ is a locally Hilbert space and
$\cal L (H_{\ell}) =\varprojlim_{p\in
Prj_{cro}(A_a)} (B(H_p), \bar\pi_{pq}^a)$ is a \prow ,
\item The morphisms $\pi_{pq}^a:A_q^a\to A_p^a$
are the restrictions of the morphisms
$\bar\pi_{pq}^a :B(H_q^a)\to B(H_p^a)$ and
$A_{\ell}$ is the sub-\prow of $\cal L (H_{\ell})$,
\item The mapping $x\mapsto (pxp)_{p\in Prj_{cro}(A_a)}$
is an embedding of
$A_a$ into $A_{\ell}$,
\item $A_{00}\subseteq A_0\subseteq A_b\subseteq A_a\subseteq A_\ell$ .
\qed
\end{enumerate}
\end{prop}

A natural question is that when the local and quasi local algebras
coincide. Here is one answer.

\begin{prop} Let $A_0$ be a \calg and $A_a$ be its atomic \walg . Then
TFAE:
\begin{enumerate}
\item For each $p\in Prj_{ro}(A_a)$, $z(p)\in Prj_{cro}(A_a)$,
\item $Prj_{cro}(A_a)\subseteq Prj_{ro}(A_a)$ is cofinal,
\item $A_{q\ell}\simeq A_\ell$ as \procs .
\end{enumerate}
\end{prop}
{\bf Proof} $(1)$ and $(2)$ are equivalent by the well known facts
that $p\leq z(p)$ and that $p\leq q$ implies $z(p)\leq z(q)$. The
equivalence of $(2)$ and $(3)$ follows from Lemma 2.4
.\qed

In next section we characterize those
\calgs which satisfy these equivalent
conditions (see Theorem 5.1).

We have used the relatively compact open projections as the
"neighborhoods " of our quasi local algebra.
Now we want to make sure that we have plenty of
these \pjs in $A_a$. More precisely, we want to check the following statements
about the \pjs of $A_a$:
\begin{enumerate}
\item Every open \pj majorizes a relatively compact open \pj ,
\item Every open \pj is the sum of finite rank projections,
\item $\bigvee Prj_{ro}(A_a)=1$.
\end{enumerate}

We start with the following lemma which is an immediate consequence of
[Ak70, I.3].

\begin{lemm} $Prj_{rc}(A_a)=\{p\in Prj(A_a): \exists a\in A^+\ \ ap=p\}$.\qed
\end{lemm}

The following result has been communicated to us by Professor Charles
Akemann for which we are grateful to him.

\begin{prop} Each $p\in Prj(A_a)$ is the sum of finite rank
projections in $A_a$. In particular, there exists some $q\in
Prj_{cp}(A_a)$ such that $q\leq p$.
\end{prop}
{\bf Proof} [Ake].\qed

\begin{cor} $\bigvee\{p: p\in  Prj_{cp}(A_a)\}=1$.
\end{cor}
{\bf Proof} Let $q$ be the LHS projection. If $1-q\neq 0$ then
there is a relatively \cpt \pj $p\leq 1-q$. But $p\leq \bar p \leq
q$, by construction, so $p\leq q(1-q)=0$, a contradiction.\qed

\section{Local algebras and multipliers}

Let $A_0$ be a \calg and $A_{00}=K(A_0)$ be its Pedersen ideal. Let
$A_a=zA_0 ^{**}$ be the atomic \walg which we associated to $A_0$ in
last section. Also there
we showed how to associate \prows $A_{\ell}$ and $A_{q\ell}$
to the atomic \walg $A_a$ of $A_0$. In the definition of these \prows
we used the family of open relatively \cpt (central) \pjs. In this
section we investigate the relation between these type of \pjs and the
Pedersen's ideal of $A_0$.

Let's introduce some notations. We denote by $Her(A_0), Idl(A_0)$, and
$Idl^\ell (A_0)$, the collection of all hereditary $C^*$-subalgebras,
all closed two sided ideal, and all closed left ideal of $A_0$,
respectively. We add index {\it cp} in each case to indicate that we
are restricting to those which are contained in $A_{00}$, so for
instance $Her_{cp}(A_0)$ is the collection of all hereditary
$C^*$-subalgebras of $A_0$ contained in $A_{00}$. Also for each
$J\subseteq A_0$ let's put
\begin{gather*}
span(J)=\text{linear span of} J,\\
conv(J)=\text{convex hull of} J,\\
her(J)=\{a\in A_0 ^+: a\leq b,\ \ \text{some} b\in J\},\\
sym(J)= cnv \cup_{t\in \Bbb T} tJ ,\\
idl_\ell(J)=\{a\in A_0: a^*a\in J\},\\
idl_r(J)=\{a\in A_0: aa^*\in J\}.
\end{gather*}
We say that $J$ is a {\it face} if it is hereditary (i.e. $her(J)=J$)
and convex (i.e. $conv(J)=J$). In this case $L=idl_\ell(J)$ is a left
ideal of $A_0$, $(L^* L)^+=J$, and $L^* L=span(J)$.

The relation between the lattice $Prj(A_a)$ and the dense ideal
$K(A_0)$ could be seen even by looking at the very definition of
$K(A_0)$. Indeed recall that
$$
K_0(A_0)=\{a\in A^+: ab=a, \ \text{some}\ b\in A^+\}\supseteq
\{a\in A^+:[a]\in Prj_{rc}(A_a)\}
$$
and $K(A_0)=span(her(conv(K_0(A))))$ [Pd66, III]. Also to each $p\in
Prj_{op}(A_a)$ one can associate a
hereditary $C^*$-subalgebra of $A_0$ defined by $her(p)=pA_a p \cap
A_0$. If moreover $p$ is relatively compact, then $her(p)\subseteq
K(A_0)$. Indeed
$$
K(A_0)= her(\cup_{p\in Prj_{ro}(A_a)} her(p))=sym(her(conv(\cup_{p\in
Prj_{ro}(A_a)} her(p)))).
$$
{\bf Remark}: In [Ped66, II] this has been stated for $A^{**}$ but
it is easy to see
that we have it for $zA^{**}$ also. ( cf. [Pd72, 3.8] and
[AP,2.6]).

We denote $idl(her(p))$ and $idl_\ell (her(p))$ simply by $idl(p)$ and
$idl_\ell (p)$, respectively. Also for $a\in
A_0$, $her(a)$ and $idl(a)$ simply denote the hereditary
$C^*$-subalgebra and closed two sided ideal of $A_0$ generated by $a$,
respectively.

The following lemma is contained in [Ped66, II]. It also could be
proved directly as follows.

\begin{lemm} If an open \pj $p$ of $A_a$ is relatively
\cpt then $her(p)\subseteq K_0(A_{0})\subseteq K(A_{0})$.
\end{lemm}
{\bf Proof}  By definition, there is $a\in A_0 ^+$ such that
$a\bar p=\bar p$. But $p\leq \bar p$, so $\bar p p=p$, hence
$ap=p$. Now if $b\in her(p)=pA_a p\cap A_0$, then $b=pb=apb=ab$,
so $b\in K_0(A_0)\subseteq K(A_0)$). \qed

Unfortunately the converse is not true, as on can see from the
following example communicated to us by Professor Gert K. Pedersen,
for which we are grateful to him. This example also shows that a (left)
ideal could be contained in $K(A_0)$ but yet the open \pj which
supports it not to be \cpt.

\begin{ex} Let $A_0$ be the universal \calg generated
by two projections p and q. Namely $A_0$ is the algebra of
continuous functions $f$ on $[0,1]$ taking values in the $2\times
2$ matrices, for which there are $s,t\in \Bbb C$ such that $$
f(0)=\left(\matrix t&0\\ 0&0\endmatrix\right)\ \ \ \  ,\ \ \ \
f(1)=\left(\matrix s&0\\ 0&t\endmatrix\right). $$ The two
generators correspond to the functions $$ p(t) =\left(\matrix
1&0\\ 0&0\endmatrix\right)\ \ \ \ ,\ \ \ \ q(t) =\left(\matrix
1-t&(t-t^2)^{\frac 1 2}\\ (t-t^2)^{\frac 1 2}&t\endmatrix\right).
$$ The point is that $A_0$ is not unital, but $K(A_0) = A_0$. Thus
the unit 1 in $\tilde A_0$ is not compact, although it supports a
left ideal (viz $A_0$) inside $K(A_0)=A_0$. Also $her(1)=A_0=
K(A_0)$. If we consider the left ideal generated by $p+q$ we
still get $A_0=K(A_0)$ so there is no hope for singly generated
ideals also. Note that $A_0$ is a subalgebra of $C([0,1],\Bbb
M_2(\Bbb C))=C([0,1])\bigotimes \Bbb M_2(\Bbb C)$ which is
separable.
\end{ex}

\begin{prop} There is a one to one correspondence
between the elements of each of the following group of sets:
\begin{enumerate}
\item $Prj_{op}(A_a), Her(A_0)$, and $Idl^\ell (A_0)$,
\item $Prj_{co}(A_a)$ and  $Idl(A_0)$,
\end{enumerate}
where in both cases the correspondence between the first two sets is
through a map of the form $p\mapsto her(p)$ and in the first
case the correspondence between the second and the third set is
through a map of the form $H\mapsto idl^\ell(H)$. Moreover under the
above correspondence we have

$(3)$ \ $Prj_{ro}(A_a)\subseteq Idl_{cp}^\ell (A_0)\subseteq Her_{cp}(A_0)$,

$(4)$ \ $Prj_{cro}(A_a)\subseteq Idl_{cp}(A_0)$.

\end{prop}
{\bf Proof} $(1)$ is proved in [Ped66, II,1.1] for $A^{**}$ and
follows for $A_a=zA^{**}$ as $z$ is central. To see $(2)$ we only
need to observe that $idl_\ell (p)$ is a two sided ideal iff $p$
is central. Now $(4)$ and the first inclusion of $(3)$ follow from
Lemma 5.1. To show the first inclusion of $(3)$, let
$her(p)\subseteq K(A_0)$, then, by [Cm70, 1.3.2],
$idl_\ell(her(p))=her(p).A_0\subseteq K(A_0).A_0\subseteq K(A_0)$.
\qed

Unfortunately the inclusions in part $(3)$ and $(4)$ could be strict
(see above example). Here we consider a class of ($\sigma$-unital)
\calgs $A_0$ for which this
holds. These \calgs are already studied by Huaxin Lin [Lin]. There he
calls them {\it pseudo commutative}, but we think that {\it pseudo
unital} is a better name.

\begin{defi}{\bf (Lin)} Let $A_0$ be a \calg . We say that $A_0$ is
pseudo unital if it has a countable approximate identity $(e_n)$ such
that
$$
\forall n\ \exists N>n\ \forall a\in A_0\ \ e_n a=e_n ae_N.
$$
\end{defi}

In particular $A_0$ would be $\sigma$-unital. Examples are $C_0(X)$
and $C_0(X)\bigotimes A$, where $X$ is a $\sigma$-compact \hau \ts ,
and $A$ is a unital \calg .

Also, following [Lin] we define a {\it support algebra} of a
$\sigma$-unital \calg $A_0$ to be a dense subalgebra of the form $\bigcup
(p_n A_0^{**} p_n \cap A_0)$, where $p_n$'s are \pjs in $A_0^{**}$ with
$e_n\leq p_n\leq e_{n+1}\ \ (n\geq 1)$, for some (countable)
approximate identity $(e_n)$ of $A_0$.

Now we are ready to prove the main result of this section.

\begin{theo} Let $A_0$ be a $\sigma$-unital \calg .
Consider the following statements: \begin{enumerate}
\item $A_0$ is pseudo unital,
\item $K(A_0)$ is a (the only) support algebra of $A_0$,
\item $\bigcup_{a\in K(A_0)^+}\  idl(a)=K(A_0)$,
\item For each hereditary $C^*$-subalgebra $b$ of $A_0$ contained in
$K(A_0)$, the closed two sided ideal generated by $B$ is also
contained in $K(A_0)$,
\item For each $a\in K(A_0)$, $idl(a)\subseteq K(A_0)$,
\item For each \pj $p\in Prj_{ro}(A_a)$, $idl(p)\subseteq K(A_0)$.
\item For each $p\in Prj_{cro}(A_a)$, $z(p)\in Prj_{ro}(A_a)$,
\item $\Gamma(K(A_0))\neq
\Delta(K(A_0))=M(A_0)$ and the spectrum $\hat A_0$ of $A_0$ is not
compact.
\end{enumerate}
Then the following implications hold
$$
(4)\implies (5)\implies (3)\iff (2)\iff (1)\implies (8),\ (5)\iff (6),
\ \text{and}\ (7)\implies (6).
$$
\end{theo}
{\bf Proof} The fact that $(1)\iff (2)\implies (8)$ is already
proved in [Lin]. To show that $(2)$ and $(3)$ are equivalent,
first note that the LHS of $(3)$ is a support algebra of $A_0$.
Indeed if $(e_n)$ is any countable approximate identity of $A_0$
contained in $K(A_0)$, then there are open central projections
$p_n$ in $A_0^{**}$ such that $idl(e_n)=p_nA_0^{**}\cap A_0$. Then
$e_n\leq p_n$ and $\cup_n idl(e_n)$ is a support algebra of $A_0$.
This is a dense ideal of $A_0$, so it contains $K(A_0)$, in
particular each $a\in K(A_0)^+$ is contained in some $idl(e_n)$,
which implies that $idl(a)\subseteq idl(e_n)$, and hence
$\bigcup_{a\in K(A_0)^+}\  idl(a)=\bigcup_{n}\ idl(e_n)$ is a
support algebra of $A_0$. Now the equivalence follows from the
fact that the two statements involved in $(2)$ are already
equivalent [Lin]. Now if $(4)$ holds and $a\in K(A_0)$, then
$her(a)\subseteq K(A_0)$ so, by assumption,
$idl(a)=idl(her(a))\subseteq K(A_0)$, which is $(5)$. To see that
$(5)$ implies $(6)$, let $p\in Prj_{ro}(A_a)$, then $p\leq a$, for
some $a\in K(A_0)^+$, so $idl(p)\subseteq idl(a)\subseteq K(A_0)$.

To finish the proof we show that
$(6)\implies(5)\implies (3)$.
Assuming $(6)$, let $a\in K(A_0)$, then the range \pj $p=[a]$ of
$a$ is a \cpt \pj in $A_0^{**}$, and $a\leq p$. Now recall that
$A_a=zA_0^{**}$, for a central \pj $z\in A_0^{**}$ and that $A_0$ could
be identified with $zA_0\subseteq A_a$. Therefore, replacing $p$ with
$zp$, we may assume that $p\in A_a$. Hence $idl(a)\subseteq
idl(p)\subseteq K(A_0)$, and $(5)$ follows. Next $(3)$ follows from
$(5)$ and the fact that $K(A_0)$ is contained in any dense ideal of
$A_0$. \qed

\section{Finite weight projections in (quasi) local algebras}

One of the most fascinating features of topological spaces is the
interrelation between topology and measure theory (of the
corresponding Borel $\sigma$-field). Recall that for a Radon (regular
positive Borel) measure on a topological space (like Haar measure on
\lcpt topological groups), we have that open sets
and compact sets
have non zero and finite measure, respectively. Also the whole space
(in the case of a topological group) has infinite measure, unless it
is compact. In non commutative topology, the {\it weights} are
analogues of Radon measures. Finite weight projections of \walgs are
studies, among others by H. Halpern, V. Kaftal, and L. Zoido [HKZ]. 
Here we consider the problem in the setting of (quasi) local
algebras.

Let's start with a brief introduction to the weight theory.
Let $M$ be a \walg and $\phi$ be a faithful normal semifinite (fns)
weight on $M$. Let $\cal N_\phi =\{x\in M: \phi(x^* x)<\infty\} $ and
$\cal M_\phi =span\cal N_\phi ^*\cal N_\phi$. Then $\cal N_\phi$ is a
left ideal of $M$ and $\cal M_\phi$ is weakly dense in $M$.

Now consider the situation where there is a \calg
$A_0$ such that $M=A_a=zA_0^{**}$, and $\phi$ is the extension of a
trace $\phi_0$ on $A_0$ to a fns
tracial weight on $A_a$
(as
before, $A_0$, being identified with
$zA_0$, is a $C^*$-subalgebra of $A_a$). We want to show that (relatively)
compact projections in $A_a$ have finite weight. We start with a
lemma which is of independent interest.

The following lemma is a classical result. The proof presented here is
communicated to us by Professor C. Akemann.

\begin{lemm} For each $p\in Prj_{cp}(A_a)$ there exists
$q\in Prj_{ro}(A_a)$ such that $p\leq q$.
\end{lemm}
{\bf Proof} Let $a\in A_0$ be such that $p\leq a$. Take $f\in
C[0,1]$ such that $0\leq f \leq 1$, $f=0$ on $[0,{\frac 1 2}]$,
and $f(1)=1$. Put $b=f(a)$ and let $q=[b]$ be the range projection
of $b$, then $q$ is open, $p\leq q$ and $\bar q$ is \cpt .\qed

\begin{lemm} If $a\in A$ and $ap=a$ for some $p\in
Prj_{cp}(A_a)$ then $a\in A_{00}=K(A_0)$.
\end{lemm}
{\bf Proof} We may assume that $a\in A_0^+$. Since $p$ is \cpt ,
there is $b\in A_0^+$ such that $bp=p$. Hence
$ba=b(ap)=b(pa)=(bp)a=pa=a$, hence $a\in K(A_0)$.\qed

\begin{lemm} If $p\in Prj_{cp}(A_a)$ then there is
$a\in K(A_{0})$ with $ap=p$.
\end{lemm}
{\bf Proof} Let $q$ be as in Lemma 6.1. By Akemann-Urysohn
lemma [Ake70, thm I.1], there is an $a\in A_0^+$ such that
$a(1-q)=0$ and $ap=p$. Now $a\bar q\geq aq=a$ and $\bar q\leq 1$,
so $a\bar q=a$ and Lemma 6.2 applies. (An alternative
proof goes as: Let $b$ be as in Lemma 6.2, then $b\in
K(A_0)$ and $p\leq b$ ).\qed

Next result is an analogue of the classical fact that \cpt sets have
finite Borel measure.

\begin{cor} With the above notation, \begin{enumerate}
\item $Prj_{cp}(A_a)\subseteq \cal M_\phi$ ,
\item $A_p^a \subseteq \cal M_\phi\ \ (p\in Prj_{cp}(A_a))$.
\end{enumerate}
\end{cor}
{\bf Proof} Let $p\in Prj_{cp}(A_a)$. By Lemma 6.3 there
is $a\in A_{00}$ such that $ap=p$. Then since $A_{00}\subseteq
\cal M_{\phi_0}$ [Ped66, III], hence $aa^*\in A_{00}\subseteq \cal
M_{\phi_0}\subseteq \cal M_\phi$, and so $a^*\in \cal N_\phi$. But
this is a left ideal of $A_a$ and we get $p=ap=pa^*\in \cal
N_\phi$. Now $p^*p=p$, so $p\in \cal M_\phi$. Also, using the fact
that $\cal N_\phi$ is a left ideal of $A_a$ again, we get
$A_p^a=pA_a p\subseteq \cal N_\phi$, from which $(2)$ implies.\qed

We don't know if the above result holds when $\phi$ is not a
trace. The only place which needs the tracial property is the
inclusion $A_{00}\subseteq \cal M_{\phi_0}$. It is known that the
canonical weight on a group algebra $C^*(G)$, for $G$ a discrete
infinite group, obtained by evaluation at the identity, fails to have
this property [Pdr]. (Note that the Haar weight of $C_r ^*(G)$ is a
trace in this case.)

We show this is the case, however, under some conditions,
for the case that
$A_a$ is either {\it finite} or {\it semifinite}
and {\it properly infinite} (see [Dix] for definitions). For
this purpose we use a deep result of Halpern, Kaftal, and Zoido 
[HKZ,thm1] which asserts that if $M$ is a properly infinite and semifinite
\walg with no type I direct summand and $Prj(\cal M_\phi)$
is a {\it lattice}, then there is a central
projection $e\in M$ such that $M_{1-e}$ is finite and $Prj(\cal
M_\phi(e.))=Prj(\cal M_\tau)$, where $\tau$ is the canonical
trace.

\begin{prop}
Let $A_0$ be a \walg and $\phi$ be a faithful normal semifinite
weight on $A_a$. Assume that either $A_a$ is finite or properly infinite,
semifinite. Also in the second case we assume that $Prj(\cal M_\phi)$
is a lattice. Then
\begin{enumerate}
\item $Prj_{cp}(A_a)\subseteq \cal M_\phi$ ,
\item $A_p^a \subseteq \cal M_\phi\ \ (p\in Prj_{cp}(A_a))$.
\end{enumerate}
\end{prop}
{\bf Proof} This follows from Corollary 6.1, and above mentioned [HKZ, thm 1].\qed

One of the reasons that we used \pros was to "expand"
the original \calg so that we can accommodate some unbounded
elements. In the theory of \walgs we usually come across such 
objects, which are mostly gathered in the set $M^\eta$ of {\it
affiliated } elements to $M$. But $M^\eta$ is merely a set (with no
apparent algebraic structure). It is desirable to have some algebras
which contain interesting unbounded operators related to a \walg . The
(quasi) local algebras are one candidate. One of the common cases where
unbounded elements affiliated to a \walg come into the play is the
Radon-Nykodym derivative of weights.
Let $M$ be a \walg and $\phi,\psi$ be fns weights on
$M$. Then there is a {\it cocycle} $u_t$ in $M$ with
\begin{gather*}
u_{t+s}=u_t \sigma_t ^\phi (u_s) \\
\sigma_t ^\psi (x)=u_t\sigma_t ^\phi (x)u_t
\end{gather*}
for each $s,t\in \Bbb R$ and $x\in M$ [Con]. (See next section for more
details about $\sigma_t ^\phi$). Moreover, if $\phi=\phi\circ\sigma_t
^\phi \ \ (t\in \Bbb R)$, then this cocycle is a one parameter group
of unitaries in the {\it fixed algebra} $M^\phi =\{x\in M:\sigma_t
 ^\phi (x)=x\ (t\in \Bbb R)\}$. Then Stone's theorem will provide us a
positive self-adjoint operator $h$ affiliated with $M^\phi$ such that
$u_t=h^{it} \ \ (t\in \Bbb R)$. Put $h_\varepsilon=h(1+\varepsilon
h)^{-1}$, then $h_\varepsilon\in M^\phi$ and
$$
\psi(x)=\lim_{\varepsilon\to 0} \phi(h_\varepsilon ^{\frac 1
2}xh_\varepsilon ^{\frac 1 2})\ \ (x\in M^+)
$$
We write $\psi =\phi(h.)$ and call $h$ the {\it Radon-Nykodym
derivative} of $\psi$ with respect to $\phi$ [PT] and
is denoted by $\frac{d\psi}{d\phi}$. If
$\psi\leq \phi$ then $0\leq h\leq 1$. If $\sigma_t ^\phi= \sigma_t
^\psi$ for all $t\in \Bbb R$, then $h$ is affiliated with $Z(M)$.

Now let $M=A_a$ be the atomic \walg of a \calg $A_0$. In this case,
we show that $h$ belongs to the quasi algebra of $A_0^{\phi_0}$. For $p\in
Prj_{ro}(A_0)$, $\phi_p(pxp)=\phi(x)$
is a faithful, semifinite, normal weights on
$A_p^a$. To avoid any
confusion we need the following lemma. The proof is similar to Lemma 2.13 of [Am] and
is omitted.

\begin{lemm} $Prj_{ro}(A_a)=Prj_{ro}(A_a^\phi)$,
and $(A_p^a)^{\phi_p}=(A_a^\phi)_p\ \
(p\in Prj(A_a))$ .\qed
\end{lemm}

Note that $\pi_{pq}^a((A_q^a) ^\phi)\subseteq (A_p^a) ^\phi$,
hence we can form the corresponding inverse limit
$A_{q\ell} ^\phi =\varprojlim_{p} (A_a\phi)_p$.

\begin{prop} Let $\phi,\psi$ be faithful, semifinite, normal
weights on $A_a$ such that $\phi=\phi\circ\sigma_t ^\phi \ \ (t\in
\Bbb R)$. Then $\frac{d\psi}{d\phi}\in A_{q\ell}^\phi \subseteq
A_{q\ell}$.
\end{prop}
{\bf Proof} Let $h=\frac{d\psi}{d\phi}\eta M^\phi$, then we know
from the construction that $h=xy^{-1}$ for some $x,y\in M^+$ [PT,
??]. But, for each $\varepsilon >0$ and $p\in Prj_{cc}(M)$,
$(\varepsilon +py)^{-1}\in M_p$, so as $\varepsilon\to 0$ we get
$(py)^{-1}\in (M_p^\phi)^\eta =M_p^\phi$, by [Am , lem 2.13]. Hence $y^{-1}\in
\check M_c ^\phi$ and so does $h$.\qed

\section{Local structure of Kac algebras}

Let $M$ be a \walg and $M \bar{\bigotimes} M$ be the {\it von Neumann
algebra tensor product} of $M$ with itself. Let $\delta:M\to
M\bar{\bigotimes} M$ be a {\it normal} injective unital $*$-homomorphism
which is {\it co-associative}, that's $(id\otimes\delta)\delta
=(\delta\otimes id)\delta$. Let $\kappa:M\to M$ be a unital involutive
anti-automorphism such that
$\delta\kappa=\varsigma(\kappa\otimes\kappa)\delta$, where $\varsigma
: M\bar{\bigotimes} M \to M\bar{\bigotimes} M $ is
the canonical {\it flip}.
Finally let $\phi:M^+\to [0,\infty]$ be
a  faithful, semifinite, normal weight on $M$. Put
$$
\cal N=\cal N_\phi =\{x\in M: \phi(x^* x)<\infty\}.
$$
This is a left ideal of $M$, which is also a {\it pre-Hilbert space}
under the inner product
$$
<x,y>=\phi(y^* x)\ \ \ (x,y\in \cal N_\phi ).
$$
Let $H_\phi$ be the Hilbert space completion of $\cal N_\phi $. We
identify $\cal N_\phi $ with a subspace of $H_\phi$. Then $\cal A_\phi
=\cal N_\phi\cap \cal N_\phi ^*$ is dense in $H_\phi$. $\cal A_\phi$
is a {\it left Hilbert algebra} under the multiplication and
involution of $M$. Also the involution of $\cal A_\phi$ is a {\it
preclosed} mapping in $H_\phi$. Let $S_\phi\in C(H_\phi)$ be its
closure and consider the {\it polar decomposition} $S_\phi
=J_\phi\Delta_\phi ^{\frac{1}{2}}$, where $J_\phi:H_\phi\to H_\phi$ is a
anti-linear isomorphism, and $\Delta_\phi$ is a positive operator,
called the {\it modular operator}. Moreover, for all $t\in \Bbb R$, we
have $\Delta_\phi ^{it}M\Delta_\phi ^{-it}=M$, and so we get a one
parameter group of automorphisms of $M$ which is called the {\it modular
automorphism group} defined by
$$
\sigma_t ^{\phi}(x)=\Delta_\phi ^{it}x\Delta_\phi ^{-it}\ \ \ (x\in
M,t\in \Bbb R)
$$
We have
$$
\phi=\phi\circ \sigma_t ^{\phi}.
$$
Now $id\otimes\phi$ is a vector valued weight, which could be
easily checked that it is faithful, semifinite,, and normal. Let $\cal
N_{id\otimes\phi}$ be the corresponding left ideal of
$M\bar{\bigotimes} M$, then we say that $\phi$ is a {\it left Haar
weight} if it satisfies the following conditions
\begin{gather*}
\delta(\cal N_\phi)\subseteq \cal N_{id\otimes\phi} \\
(id\otimes\phi)((1\otimes y^*
)\delta(x))=\kappa((id\otimes\phi)(\delta(y^* )(1\otimes x)))\ \ \
(x,y\in \cal N_\phi) \\
\kappa\sigma_t ^{\phi}=\sigma_{-t} ^{\phi} \kappa .
\end{gather*}
Then $\phi$ is {\it left invariant} in the sense that
$$
(id\otimes\phi)\delta(x)=\phi(x)1\ \ \ (x\in M^+ ).
$$

Under the above conditions $K=(M,\delta,\kappa,\phi)$ is called a {\it
Kac algebra}. Also $\delta,\kappa,\phi$ are called the {\it co-product,
co-inverse} and {\it Haar measure } of $K$, respectively.

Two classical examples of Kac algebras are the {\it commutative } and
{\it co-commutative} \kacs $K_a=(L^\infty (G),\delta_a ,\kappa_a
,\phi_a )$ and $K_s=(\cal L (G),\delta_s ,\kappa_s ,\phi_s )$, where $G$
is a \lcpt topological group, $L^\infty (G)$ is the algebra of all
{\it essentially bounded} complex \fns on $G$, and
$$
\cal L(G)=\{\lambda(f): f\in L^1 (G)\}^{''}
$$
 is the {\it (left) group von Neumann algebra} of $G$, and the other
ingredients are defined by
\begin{gather*}
\delta_a (f)(s,t)=f(st),\ \ \delta_s (\lambda(g))=
\lambda(f)\otimes\lambda(f)\\
\kappa_a(f)(s)=\tilde f (s)=f(s^{-1}),\ \
\kappa_s(\lambda(g))=\lambda(\tilde g)^* \\
\phi_a(f)=\int_{G} f(s)ds,\ \ \phi_s(h)=h(e)
\end{gather*}
where $f\in L^\infty (G), g\in L^1 (G)$, and $h\in L^1(G)\cap L^2
(G)$.

In this section we want to introduce the analogs of {\it local
neighborhoods} for \kacs . In the previous sections we observed that
(open) projections can play this role for us. One should note however
that these concepts were defined in the context of {\it atomic} \walgs
. It is possible to associate an atomic \walg to each \kac $K$ (it
would be a quotient of $W^*(K)$ in the notation of [ES]), but there is
a more suitable (and natural) way to handle this. Starting with a
\calg $A_0$, one main reason that the non commutative topology works
is the fact that there is an isometrically isomorphic copy of $A_0$ in
its atomic \walg $zA_0^{**}$ [Pd72, 3.8]. Now we have a similar
situation in \kacs , namely given a \kac \K , let $A_0=C_r ^*(\hat K)$
be the \calg completion of $\hat\lambda(\hat M_*)$ [ES]. Then
$A_0\subset M$ has the same norm. We define the concepts of non
commutative topology in this framework and show that we have plenty of
open, closed, and \cpt \pjs inside $M$.

In what follows we freely use the notations and terminology of
[ES]. Let \K be a \kac and $\hat K$ be its dual \kac . Let $\hat
\lambda : \hat M_*\to M\subseteq B(\cal H_{\hat\phi})$ be the
Fourier transform of $\hat K$ and $A_0=C_r ^*(\hat K)$ be the
\calg completion of $\hat\lambda(\hat M_*)$. Then $A_0\subseteq
M=\hat\lambda(\hat M_*)^{''}$, where the commutants are with
respect to the Fourier \repn . On the other hand, the inclusion
map $j:A_0 \to M$ extends to a normal \shomo $j ^{''}: A_0
^{''}\to M^{''}$, where the commutants this time are with respect
to the universal \repns . But then $A_0 ^{''}$ is the same as the
second dual $A_0 ^{**}$ [Sak, 1.17.2], and $M^{''}=M$. Let
$I=ker(j ^{''})$, then $I$ is a $weak^*$-closed two sided ideal of
$A_0 ^{**}$, so there is a central \pj $w\in A_0 ^{**}$, such that
$I=(1-w) A_0 ^{**}$.

\begin{lemm} With the above notations, $M\backsimeq wA_0 ^{**}$, as
\walgs . Indeed the restriction of $j^{''}$ to $wA_0 ^{**}$ is
injective and so an isomorphism. Also the following diagram
commutes.

$$\begin{matrix}
A_0 & \overset{j}{\lo} & \omega A_0 \\
\mapdown{i} && \mapdown{} \\
M & \overset{\simeq}{\lo} & \omega A_0^{**}
\end{matrix}$$

where the top map is the \shomo which sends $a\in A_0$ to $wa$ and
$i$ and $j$ are the inclusion maps.
\end{lemm}
{\bf Proof} The first statement follows from the fact that $j
^{''}$ is surjective (since $M=\hat\lambda(\hat M_*)^{''}$) and
that $ker(j ^{''})=(1-w) A_0 ^{**}$. To see that the diagram
commutes, take any $a\in A_0$, then this is sent to $wa\in wA_0
^{**}$ in the top route . Now $j^{''}$ is a unital \homo so
$j^{''}(1)=w$. Also $j^{''}(a)=j(a)=a$, therefore the isomorphism
on the bottom (which is the inverse of the restriction of $j^{''}$
to $wA_0 ^{**}$) sends $a$ to $wa$ and the diagram commutes. \qed

Given $x\in M$ and $J\subseteq M$, following [GK], we put
$$
e(x)=\wedge \{p\in Prj(M): xp=x\}\ \ \text{and}\ \
e(J)=\vee\{e(x): x\in J\}.
$$
Then $e(x)$ and $e(J)$ are \pjs in $M$. Similarly we define central
\pjs $z(x)$ and $z(J)$ in $M$ by
$$
z(x)=\wedge \{p\in Prj_{cn}(M): xp=x\}\ \ \text{and}\ \
e(J)=\vee\{z(x): x\in J\}.
$$

\begin{defi} Let \K be a \kac and $A_0=C_r ^*(\hat K)$.
A \pj $e\in Prj(M)$ is called open if there is a set
$J\subseteq A_0$ such that $e=e(J)$. It is called closed if $1-e$ is
open. Also $e\in Prj(M)$ is called relatively \cpt if there is $a\in
A_0$ such that $0\leq a\leq 1$ and $e\leq a$. Finally a closed relatively
\cpt \pj is called \cpt . We denote the set of open, closed,
relatively \cpt ,and \cpt \pjs of $M$ by $Prj_{op}(M)$, $Prj_{c\ell}(M)$,
$Prj_{rc}(M)$, and $Prj_{cp}(M)$, respectively. We also freely use the
abbreviations of Notation~\ref{notation} in this context.
\end{defi}

Next we show that there are plenty of open (and so closed) \pjs in
$M$. The few next results are already proved in [GK] in a different
context (where $M=zA_0 ^{**}$). We give detailed proofs here, with
some slight modifications when needed, to make
sure they work in our context. In what follows \K is a \kac and $A_0
=C_0 (K)\subseteq M$.

\begin{prop} Let $a\in A_0$ be self adjoint and $U\subseteq \Bbb R$
be open. Let
$$
a=\int_{-\infty}^{\infty} sdE(s)
$$
be the spectral decomposition of $a$ and put $e=\int_{U} dE(s)$. Then
$e\in Prj_{op}(M)$. Moreover if $U$ has \cpt closure in $\Bbb R$, then
$e$ is relatively \cpt .
\end{prop}
{\bf Proof} $e$ is clearly a \pj . The fact that $e\in M$ follows
from the fact that $E(s)\in M$ for each $s\in\Bbb R$ [Sak,
1.11.3]. Take any $f\in C(\Bbb R)$ such that $f=0$ off $U$ and
$f>0$ on $U$. Put $b=f(a)$. Then $b\in \tilde A_0$, the (minimal)
unitization of $A_0$.

Let's observe that $e=e(b)$. Indeed $b=f(a)=\int_{-\infty}^{\infty}
f(s)dE(s)$ and $e(b)$ belongs to the \walg generated by $b$ [Sak,
1.10.4] and so to the one generated by $\{E(s) : s\in \Bbb R\}$. In
particular $e(b)= \int_{-\infty}^{\infty} g(s)dE(s)$, for some bounded
\fn $g$. By definition of $e(b)$ then we have $g^2=g=gf$, which means
that $g$ is the characteristic \fn of a set containing $U$, and so by
minimality, $g=\chi_U$ and $e(b)=e$.

Take any \ai $(e_\alpha)$ of $A_0$ consisting of
positive elements. Then $b_\alpha =b^{\frac  1 2}e_\alpha b^{\frac  1
2}\in A_0$, as $A_0$ is an ideal of $\tilde A_0$. Let's show that
$\vee e(b_\alpha)\leq e(b)=e$. Let $p\in Prj(M)$ and $bp=0$ , then
$(b^{\frac  1 2}p)^* (b^{\frac  1 2}p)=0$ and so $b^{\frac  1 2}p=0$
which implies $b_\alpha p=0$, for each $\alpha$, which means
$\vee e(b_\alpha)\leq e(b)$. Conversely, if $b_\alpha p=0$, for each
$\alpha$, then by the fact that $e_\alpha\to 1$
and $b_\alpha\to b$ in weak-operator topology [Sak, 1.7.4], and that
multiplication of $M$ is separately weakly \cnt , we get $bp=0$, and
so $e(b)p=0$, which means that $e(b)\leq \vee e(b_\alpha)$. Therefore
$e=e(b)=\vee e(b_\alpha)=e(\{b_\alpha\})$ is open.

Finally if $U$ is relatively \cpt in $\Bbb R$, then we can choose a
\fn $g\in C_{00}(\Bbb R)^+$ such that $0\leq g\leq 1$ and $g=1$ on $U$, then
$c=\int_{-\infty}^{\infty} g(s) dE(s)\in A_0$ satisfies $0\leq c \leq
1$ and
$$
e=\int_{-\infty}^{\infty} \chi_U (s) dE(s)\leq \int_{-\infty}^{\infty}
g(s) dE(s)=c.\qed
$$

For each $e\in Prj(M)$ let's consider
$$
L(e)=\{a\in A_0: ae=a\},
$$
which is a norm-closed left ideal of $A_0$, and note that $ae=a$ iff
$e(a)\leq e$, for each $a\in A_0$, therefore
$$
L(e)=\{a\in A_0: e(a)\leq e\}=\{a\in A_0: e(a^* a)\leq e\}.
$$
Recall that for a \pj $e$, the {\it interior} $e^o$ of $e$ is defined
by $e^o =\vee\{p\in Prj_{op}(M): p\leq e\}$. The {\it closure} $\bar
e$ of $e$ is defined similarly. Next let's observe that $L(e)=L(e^o)$,
for each $e\in Prj(M)$. Indeed, for each $a\in L(e)$, $e(a)\leq
e$, and so $e(a)\leq
e^o$, as $e(a)$ is open. Hence $L(e)\subseteq L(e^o)$. The converse
follows from $e^o \leq e$.

Next, for each $a\in L(e^o)$, $e(a)\leq
e^o$, so $e(L(e^o))\leq e^o$. Also $e_o$ is open, so $e^o=e(J)$, for
some $J\subseteq A_0$, but then $J\subseteq L(e^o)$ (by definition of
$e(J)$) and so $e^o\leq e( L(e^o))$. Therefore
$$
e( L(e))=e( L(e^o))=e^o .
$$

\begin{lemm} With the above notation, TFAE:
\begin{enumerate}
\item For every closed left ideal $L$ of $A_0$, $e(L)=1$ implies $L=A_0$,
\item For every maximal closed left ideal $J$ of $A_0$, $L(e(J))=J$.
\end{enumerate}
\end{lemm}
{\bf Proof} If $(1)$ holds and $J$ is maximal, then for each $x\in
J$ we have $x\leq e(J)$ and so $x\in L(e(J))$, by definition, i.e.
$J\subseteq L(e(J))$. Hence, by maximality, either $J=L(e(J))$ or
$L(e(J))=A_0$. But the second equality implies that $e(J)=1$, and
so $J=A_0$, by $(1)$, which is a contradiction. Conversely, if
$(2)$ holds, and $L$ is a closed left ideal with $e(L)=1$ but
$L\neq A_0$, then there is a (proper) maximal ideal $J\supseteq
L$. Then $e(J)\geq e(L)=1$, and so $e(J)=1$, which implies by
$(2)$ that $J=L(e(J))=L(1)=A_0$, which is a contradiction. \qed

\begin{defi} Let \K be a \kac and $A_0=C_r ^*(\hat K)$. We say that $K$
satisfies condition (D) if any of the equivalent conditions of above
lemma holds.
\end{defi}

The above lemma is one of the departure points of our theory and the
Akemann-Giles non commutative topology. In their context the condition
$(D)$ is automatically satisfied as their \walg is atomic. For \kacs ,
however, this condition does not hold in general even for the
commutative case. If $M=L^\infty (\Bbb R)$ then $A_0=C_0(\Bbb R)$, and
we have an uncountable number of distinct maximal closed ideals of
$A_0$ whose support is $1$, namely the ideals
$$
J_x=\{f\in C_0(\Bbb R): f(x)=0\}\ \ \ (x\in \Bbb R).
$$
However, if we replace $\Bbb R$ with $\Bbb R_d$ ($\Bbb R$ with the
discrete topology) then we get $e(J_x)=1-\delta_x$ and condition $(D)$
is satisfied. This is indeed true for any discrete \kac .
Recall that a \kac \K is called {\it discrete} if the Banach algebra
$M_*$ is unital [ES].

\begin{lemm} Any discrete \kac satisfies condition (D).
\end{lemm}
{\bf Proof} If \K is discrete, then $M$ is an atomic \walg [ES,
6.3.4] and we can apply the results of [GK] to our case. In
particular, condition $(D)$ follows from the proof of [GK, 5.5].
\qed

\begin{lemm} With the above notation
\begin{enumerate}
\item If $J$ is a maximal closed left ideal of $A_0$, then $1-e(J)$ is
a minimal closed \pj (possibly zero) of $M$
(i.e. it is closed and minimal among all closed \pjs),
\item If condition $(D)$ holds, then
for each closed left ideal $L$ of $A_0$, $L=L(e(L))$.
\end{enumerate}
\end{lemm}
{\bf Proof} We start by showing half of $(ii)$. Given $x\in L\in
Idl^\ell (A_0)$, we have $x\leq e(L)$, and so $x\in L(e(L))$, i.e.
$L\subseteq L(e(L))$. Now to show $(i)$, take any $f\in
Prj_{c\ell}(M)$ such that $f\leq 1-e(J)$, then $1-f\geq e(J)$, and
so $L(1-f)\supseteq L(e(J))\supseteq J$. By maximality, we get
$L(1-f)=J$. Hence $e(J)=e(L(1-f))=(1-f)^o = 1-f$, which proves
$(i)$. To prove $(ii)$, assume that $L$ is a closed left ideal
with $L\neq L(e(L))$, then there is a maximal closed left ideal
$J$ of $A_0$ such that $J\supseteq L$ but $J\not\supseteq L(e(L))$.
Given $a\in L$ we have $x\in J$, and so $xe(J)=x$, i.e.
$x(1-e(J))=0$. This holding for each $x\in L$, we get
$e(L)(1-e(J))=0$. But then $L(e(L))(1-e(J))=\{0\}$, by definition.
Hence for each $x\in L(e(L))$, we have $x=xe(J)$ and so $x\in
L(e(J))$. But  $L(e(J))=J$, by condition $(D)$, therefore
$J\supseteq L(e(L))$, which is a contradiction. \qed

\begin{prop} The map $e\mapsto L(e)$ is an order preserving
injection of $Prj_{op}(M)$ into $Idl^\ell (A_0)$. Moreover, if $K$
satisfies condition $(D)$, this is a bijection.
\end{prop}
{\bf Proof} For $e,f\in Prj_{op}(M)$, then $e=e^o$ and $f=f^o$.
Therefore if $e\leq f$, then $e(L(e))\leq e(L(f))$, and so
$L(e)\subseteq L(f)$. Also if $L(e)=L(f)$, then $e=e( L(e))=e(
L(f))=f$. To show that this map is onto, if $K$ satisfies
condition $(D)$, we only need to observe that in this case,
$L=L(e(L))$, for each $L\in Idl^\ell (A_0)$, by above lemma. \qed

\begin{cor} The map $e\mapsto L(e)$ is an order preserving
injection of the set $Prj_{co}(M)$ of all central open \pjs of $M$
into $Idl (A_0)$. Moreover, if $K$ satisfies condition $(D)$ this
is a bijection.\qed
\end{cor}

As we pointed out earlier the \cpt \pjs in non commutative topology do
not behave in some cases like the characteristic \fns of \cpt sets
in the commutative case. For instance an Akemann's \cpt \pj $p$ could be
majorized by a union of open \pjs where non of its finite subunions
could majorize $p$ (think of a rank one \pj covered by union of rank
one \pjs ). Giles and Kummer showed that this pathology doesn't occur
in the {\it \cpt} case [GK, 5.7]. We use a simplified version of their
proof to show the
implication $(3)\implies (4)$ in the following result.

\begin{theo} Let \K be \kac . Consider the following conditions:
\begin{enumerate}
\item $K$ is \cpt , i.e. $\phi$ is finite,
\item $A_0=C_r ^*(\hat K)$ is unital,
\item $1\in Prj_{cp}(M)$,
\item (finite covering property) For each family
$(e_i)_{i\in I}$ of open \pjs of $M$ with
$\bigvee_{i\in I} e_i=1$, there is a finite subset $J\subseteq I$ such
that $\bigvee_{i\in J} e_i=1$,
\item (finite intersection property) For each family
$(e_i)_{i\in I}$ of closed \pjs of $M$ with
$\bigwedge_{\ i\in I} e_i=0$, there is a finite subset $J\subseteq I$ such
that $\bigwedge_{\ i\in J} e_i=0$.
\end{enumerate}
Then $(1)\iff (2)\iff (3)$ and $(4)\iff (5)$. Moreover, if $K$
satisfies condition $(D)$, then we have $(3)\implies (4)$.
\end{theo}
{\bf Proof} The If $K$ is \cpt then the Banach algebra $\hat M_*$
is unital [ES, 6.2.5]. But $A_0=C_r ^*(\hat K)$ is the
$C^*$-completion of $\hat M_*$, so $A_0$ is also unital.
Conversely, if $A_0$ is unital, then $K$ is \cpt by the proof of
[ES, 6.3.3]. Next we show that $(2)$ and $(3)$ are equivalent. If
$1$ is \cpt, then $1a=a$, for some $a\in A_0$, i.e. $A_0$ is
unital . Conversely, if $A_0$ is unital, then $1\in A_0$ is a \cpt
\pj of $M$. The equivalence of $(4)$ and $(5)$ is trivial.

Now assume that $K$ satisfies condition $(D)$,
and that $1\in A_0$ is \cpt. Then consider the
left ideal $L=\sum_{i\in I} L(e_i)$. Then $e(L)\geq e(L(e_i))=e_i
^o=e_i$, for each $i$, and so $e(L)\geq\vee e_i =1$. Hence $e(L)=1$
and so $L=A_0$ by condition $(D)$. In particular $1\in L=\sum_{i\in
I} L(e_i)$. Therefore there is a finite subset $J\subseteq I$ such
that $1\in L=\sum_{i\in J} L(e_i)$. Let $e_J=\bigvee_{i\in J} e_i $,
then $e_J=1e_J=1$, and $(4)$ is proved.
\qed

{\small {\bf Acknowledgement}: This paper is part of the author's
Ph.D. thesis in the University of Illinois at Urbana-Champaign
under the supervision of Professor Zhong-Jin Ruan. I would like to
thank him for his moral support and scientific guidance during my
studies.}

\end{document}